\documentclass[12pt]{amsart}
\usepackage{epsfig}

\title[Noncommutative trigonometry and the A-polynomial]{Noncommutative
 trigonometry and the A-polynomial of the trefoil knot}
\author{R{\u{a}}zvan Gelca}
\address{Department of Mathematics, 
University of Michigan, Ann Arbor, MI 48109 and Institute of Mathematics
of the Romanian Academy, Bucharest, Romania}
\email{rgelca@math.lsa.umich.edu}

\newtheorem{theorem}{Theorem}

\newtheorem{lemma}{Lemma}

\newcommand{\lcr}{\raisebox{-5pt}{\mbox{}\hspace{1pt}
                  \epsfig{file=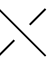}\hspace{1pt}\mbox{}}}
\newcommand{\ift}{\raisebox{-5pt}{\mbox{}\hspace{1pt}
                  \epsfig{file=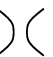}\hspace{1pt}\mbox{}}}
\newcommand{\zer}{\raisebox{-5pt}{\mbox{}\hspace{1pt}
                  \epsfig{file=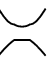}\hspace{1pt}\mbox{}}}

\begin{document}
\maketitle

\section{Introduction}

In \cite{frgelo}   a noncommutative generalization of the 
A-polynomial of a knot of  Cooper, Culler, Gillet,
Long and Shalen \cite{CCGLS} was introduced.
This generalization consists of a finitely generated
 left ideal of polynomials in the quantum plane, 
the noncommutative A-ideal, and 
was defined  based on Kauffman bracket skein modules, 
by deforming the ideal generated by the A-polynomial with respect
to a parameter. The deformation was  possible because of 
 the relationship between the skein module with the variable
$t$ of the Kauffman bracket evaluated at  $-1$ and the $SL(2,C)$-character
variety of the fundamental group, which was  explained in \cite{bullock2}.   
The purpose of this paper is to show how the 
noncommutative A-ideal can be computed for the left- and right-handed
trefoil knots. As it will be seen below, this task reduces
to trigonometric computations in the noncommutative torus.

The computation of the A-ideal 
relies heavily on understanding the image in the 
skein module of the knot complement, 
of the Kauffman bracket skein algebra of the cylinder over
the boundary torus, under the map induced by the 
natural inclusion of the boundary into the knot complement.
This amounts also to understanding the action of 
the   Kauffman bracket skein algebra of the cylinder over
the torus on the Kauffman bracket skein module of the knot
complement, induced by the inclusion map. The solution for the left-handed 
trefoil to both these problems makes
the object of  Theorem 1 bellow.
In Theorems  2 and 4 we describe the peripheral ideal of the left-handed
trefoil knot (i.e. the left ideal of the skein algebra of the cylinder
of the torus consisting of elements  which become zero in the 
skein module of the knot complement). Finally, Theorem $3$ 
lists three generators of the noncommutative A-ideal of the left-handed 
trefoil
when $t$ is not an eighth root of unity. The case $t=-1$ makes
the object of Theorem 5. In the last section, we list the 
analogous results for the right-handed trefoil. 
The proofs of the results are sketchy, since we  preferred to 
skip long routine computations, and to focus only
on main ideas. 

\section{Preliminary facts}

Throughout this paper $t$ will denote a fixed complex number. 
A {\em framed link} in an orientable manifold $M$ is a 
disjoint union of annuli. In the case where the manifold is
the cylinder over the torus, framed links will be identified with 
curves, using the convention that the annulus is parallel to the framing.
Let $\mathcal{L}$ be the set of isotopy classes of links in 
the manifold $M$, including the empty link. 
Consider the complex vector space with basis $\mathcal{L}$,
and factor it by the smallest subspace containing all expressions
of the form $\displaystyle{\lcr-t\zer-t^{-1}\ift}$
and 
$\bigcirc+t^2+t^{-2}$, where the links in each expression are
identical except in a ball in which they look like depicted.
This  quotient is denoted by $K_t(M)$ and is called the Kauffman 
bracket skein module of the manifold \cite{przytycki}. 
In the case of a cylinder over
a surface, the skein module has an algebra structure induced
by the operation of gluing one cylinder on top of the other.
The operation of gluing the cylinder over $\partial M$ to $M$ induces
a $K_t(\partial M\times I)$-left module structure  on $K_t(M)$.
We denote by $*$ the multiplication in $K_t(\partial M\times I)$
and by $\cdot $ the left action of this algebra on $K_t(M)$.

Let us discuss in more detail two structure results about
the Kauffman bracket skein algebra of the torus and the Kauffman bracket
skein module of the complement of the trefoil knot.
For this we need to introduce two families of polynomials.
The first family consist of the  Chebyshev polynomials
$T_n$, $n\geq 0$,
 defined by $T_0(x)=2$, $T_1(x)=x$, and $T_{n+1}(x)=xT_n-T_{n-1}$.
Recall that they  arises when expressing $2\cos n\alpha$ as a function of 
$2\cos \alpha$.  The
second family of polynomials $S_n$, $n\geq 0$, is  
closely related to the first,
satisfying the same reccurence relation, but with $S_0(x)=1$ and $S_1(x)=x$. 
These arise when writing $\sin{(n+1)\alpha}/\sin \alpha $ as a function of
$\cos \alpha$. They satisfy $S_n=T_n+T_{n-2}+T_{n-4}+\cdots $, where the sum
ends in a $1$ (not $2$) if $n$ is even.
Extend both polynomials by the reccurence relation to
all indices $n\in{\mathbb Z}$. Note that $T_{-n}=T_n$, while
$S_{-n}=-S_{n-2}$. 

For a link $\gamma $ in a skein module we will denote by
$\gamma ^n$ the link consisting of $n$ parallel copies of 
$\gamma$, and extend the notation to polynomials.
For a pair of integers $(p,q)$, we denote by 
$(p,q)_T$ the element $T_n((p',q'))$ of the Kauffman bracket skein module
of the cylinder over the torus $K_t({\mathbb T}^2\times I)$
where $n$ is the greatest common divisor of $p$ and $q$,
$p'=p/n$, $q'=q/n$ and $(p',q')$ is the  simple closed curve of 
slope $p'/q'$ on the torus. Define analogously $(p,q)_{JW}=S_n((p'q'))$.
The index is motivated by the fact that $S_n((p'q'))$ is the 
$(p',q')$-curve colored by the $n$-th Jones-Wenzl idempotent.

As a complex vector space, $K_t({\mathbb T}^2\times I)$ is spanned
by the elements $(p,q)_T$, $p\geq 0$. 
It was proved in \cite{frogel} that $K_t({\mathbb T}^2\times I)$
 is canonically
isomorphic to the subalgebra of the noncommutative torus \cite{connes} spanned
by noncommutative cosines. The isomorphism is a consequence of 
the {\em product-to-sum} formula for skeins in $K_t({\mathbb T}^2\times I)$:
\begin{eqnarray*}
(p,q)_T*(r,s)_T=t^{|^{pq}_{rs}|}
(p+r, q+s)_T+
t^{-|^{pq}_{rs}|}
(p-r,q-s)_T.
\end{eqnarray*}
This formula will be fundamental in all computations throughout this paper.

As proved in \cite{bullock}, the Kauffman bracket skein module
of the complement of the trefoil (whether left or right),
 is the vector space generated by $x^n$ and $x^ny$, where
$n\geq 0$, and $x$ and $y$ are the curves shown in Fig. 1.
\begin{figure}[htbp]
\centering
\leavevmode
\epsfxsize=2.2in
\epsfysize=1.2in
\epsfbox{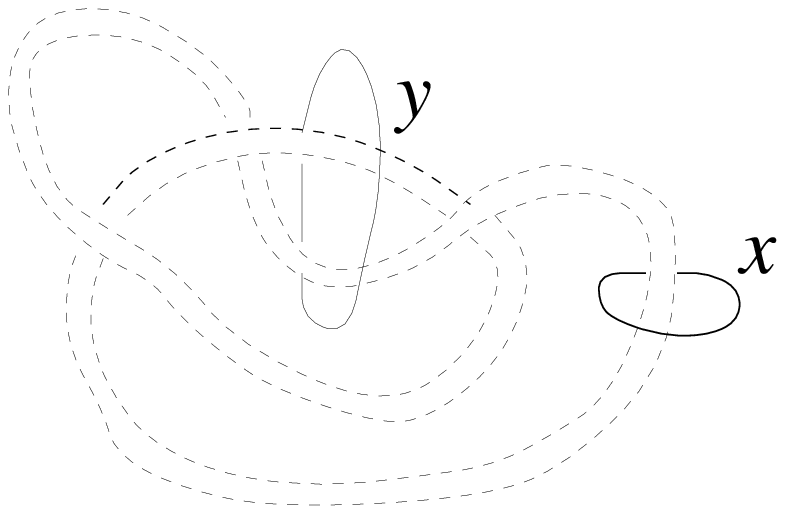}

Figure 1.  
\end{figure}
As we will see below, it is more useful to consider a different 
basis of this  vector space, namely $S_n(x)$ and $S_n(x)y$, $n\geq 0$. 
As it was in the case of the unknot \cite{frogel}, by passing from the  
boundary torus to the knot complement one
passes from curves colored by $T_n$'s to curves colored by
$S_n$'s. 

If $M$ is the complement of a knot $K$, then one can introduce a noncommutative
generalization of the A-polynomial  of $K$ defined in \cite{CCGLS},  in the
following way (see \cite{frgelo}). Denote by $\pi$ the map between 
skein modules induced
by the inclusion $\partial M \times I)\subset M$ and let $I_t(K)$ be the
kernel of $\pi$. The noncommutative A-ideal of $K$ is defined to 
be the left ideal obtained by extending $I_t(K)$ to the subalgebra
of the noncommutative torus consisting of trigonometric polynomials,
and then contracting it to the quantum plane. Recall
that the algebra of noncommutative trigonometric polynomials
is ${\mathbb{C}_t[l,l^{-1},m,m^{-1}]}$, and the quantum plane
is ${\mathbb{C}_t[l,m]}$ \cite{kassel},
  where in both,  $l$ and $m$ satisfy $lm=t^2ml$.
The classical A-polynomial is obtained by making $t=-1$, replacing
$l$ and $m$ by $-l$ and $-m$ respectively,  and taking the generator of the 
radical of the part of the ideal that has   Krull dimension equal to $1$. 
The fact that this is indeed the A-polynomial follows from Bullock's 
\cite{bullock2}
characterization of the Kauffman bracket skein module at $t=-1$ as 
the affine $SL(2,C)$-character variety ring of the fundamental group (see also
\cite{przsik}),  together with 
 the fact that $K_{-1}({\mathbb T}^2\times I)$ has
no nilpotents. 

With the above notation, the elements $e_{p,q}=t^{-pq}l^pm^q$ in
${\mathbb{C}_t[l,l^{-1},m,m^{-1}]}$ are called 
noncommutative exponentials, the elements $\cos _t(p,q)=(e_{p,q}+e_{-p,-q})/2$ 
are called noncommutative cosines, and the elements $\sin _t(p.q)=
(e_{p,q}-e_{-p,-q})/2$
are called noncommutative sines. The morphism from 
$K_t({\mathbb T}^2\times I)$ to the noncommutative torus
maps $(p,q)_T$ to $2\cos _t(p,q)$, and under this identification
the computations  become trigonometric manipulations
of  $\cos _t(p,q)$ and $\sin_t (np, nq)/\sin_t(p,q)$, hence the title
of the paper.

\section{The $K_t({\mathbb T}^2\times I)$-module structure
of the skein module of the complement of the left-handed trefoil
knot}

In this section $K$ will denote the left-handed trefoil knot
and $M$ the complement of a regular neighborhood of this knot.

\begin{lemma}
In the Kauffman bracket skein module of the complement of the
trefoil knot, the following equality holds
\begin{eqnarray*}
y^2=-t^2S_2(x)y-t^4S_2(x)+S_0(x).
\end{eqnarray*}
\end{lemma}

\proof
The element $y^2$ can be obtained by 
the Kauffman bracket skein relation for the 
unknot like in Fig. 2.  
\begin{figure}[htbp]
\centering
\leavevmode
\epsfxsize=3.7in
\epsfysize=1.2in
\epsfbox{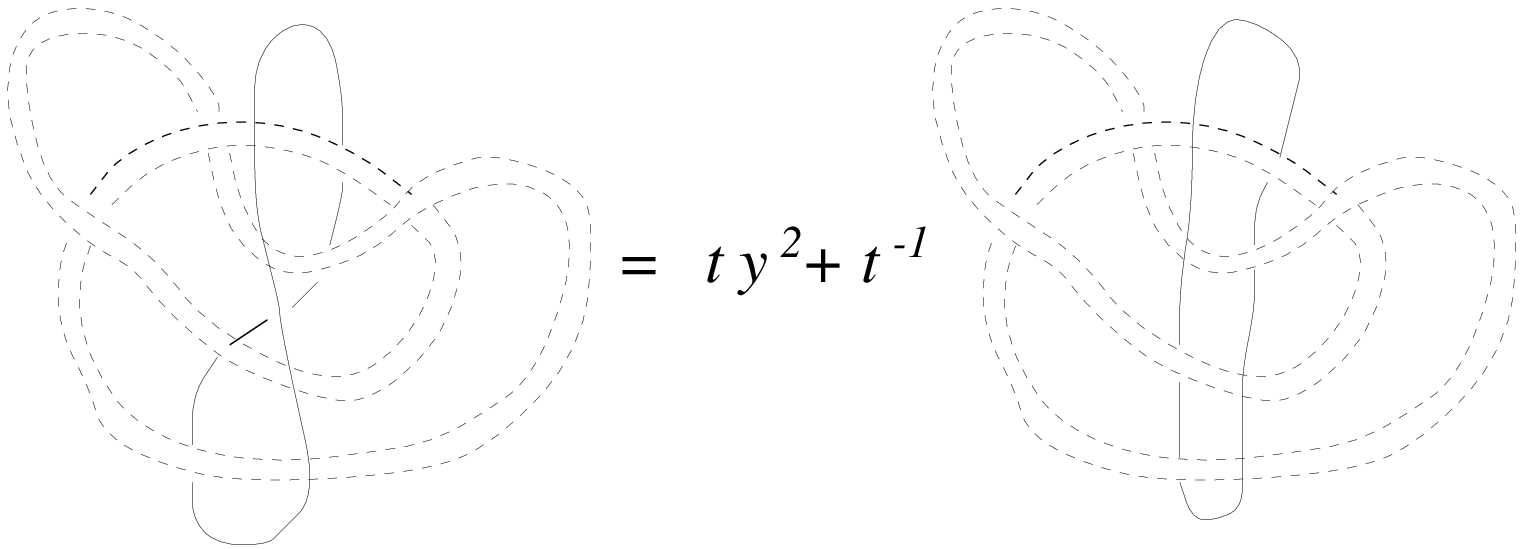}

Figure 2.  

\centering
\leavevmode
\epsfxsize=4in
\epsfysize=3.8in
\epsfbox{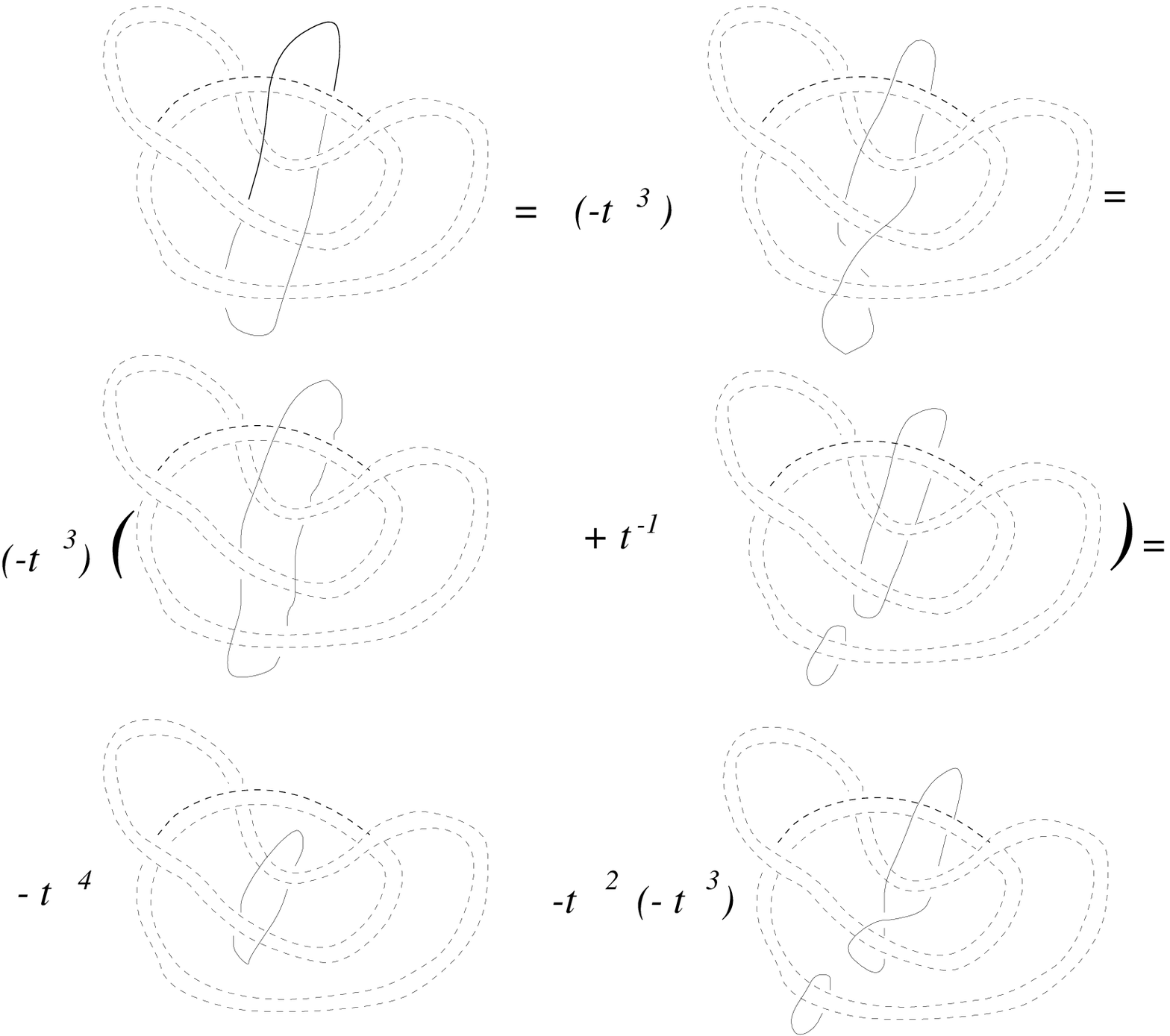}

Figure 3.  
\end{figure}
To compute the skein on the right side, one can proceed like
in Fig. 3. It is not hard to see that in the end one obtains 
$-t^4y+t^4x^2y+t^6x^2$. Since the unknot with a positive twist from
Fig. 1 is equal to $(-t^3)(-t^2-t^{-2})$, we get 
$y^2=1+t^4+t^2y-t^2x^2y-t^4x^2$. Grouping these terms appropriately we
get the formula from the statement.
\qed

\begin{lemma}
In the Kauffman bracket skein module of the complement of the
trefoil knot, the following equality holds
\begin{eqnarray*}
y^3=t^4S_4(x)y+2S_0(x)y+t^6S_4(x)+t^{10}S_0(x).
\end{eqnarray*}
\end{lemma}

\proof
To compute $y^3$ we use the relation from Fig. 4. It is not hard to
see that the skein on the left is equal to $(t+t^5)y$, while the skein on
the right can be transformed like in Fig. 5.

\begin{figure}[htbp]
\centering
\leavevmode
\epsfxsize=3.5in
\epsfysize=1.2in
\epsfbox{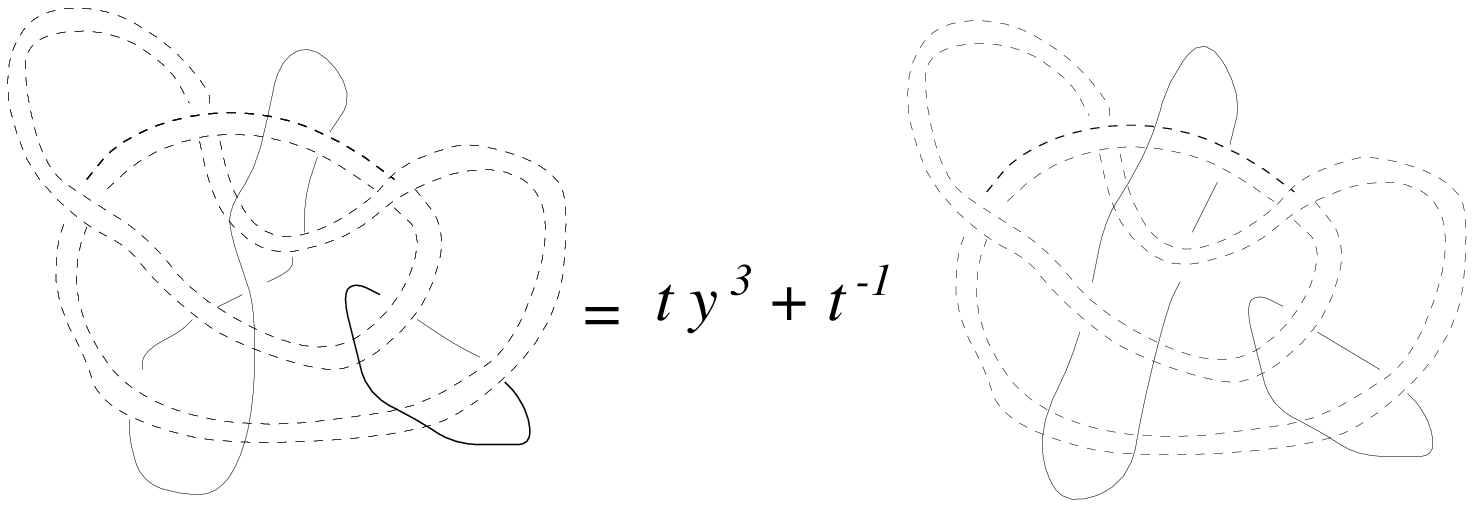}

Figure 4.  

\centering
\leavevmode
\epsfxsize=5in
\epsfysize=3.6in
\epsfbox{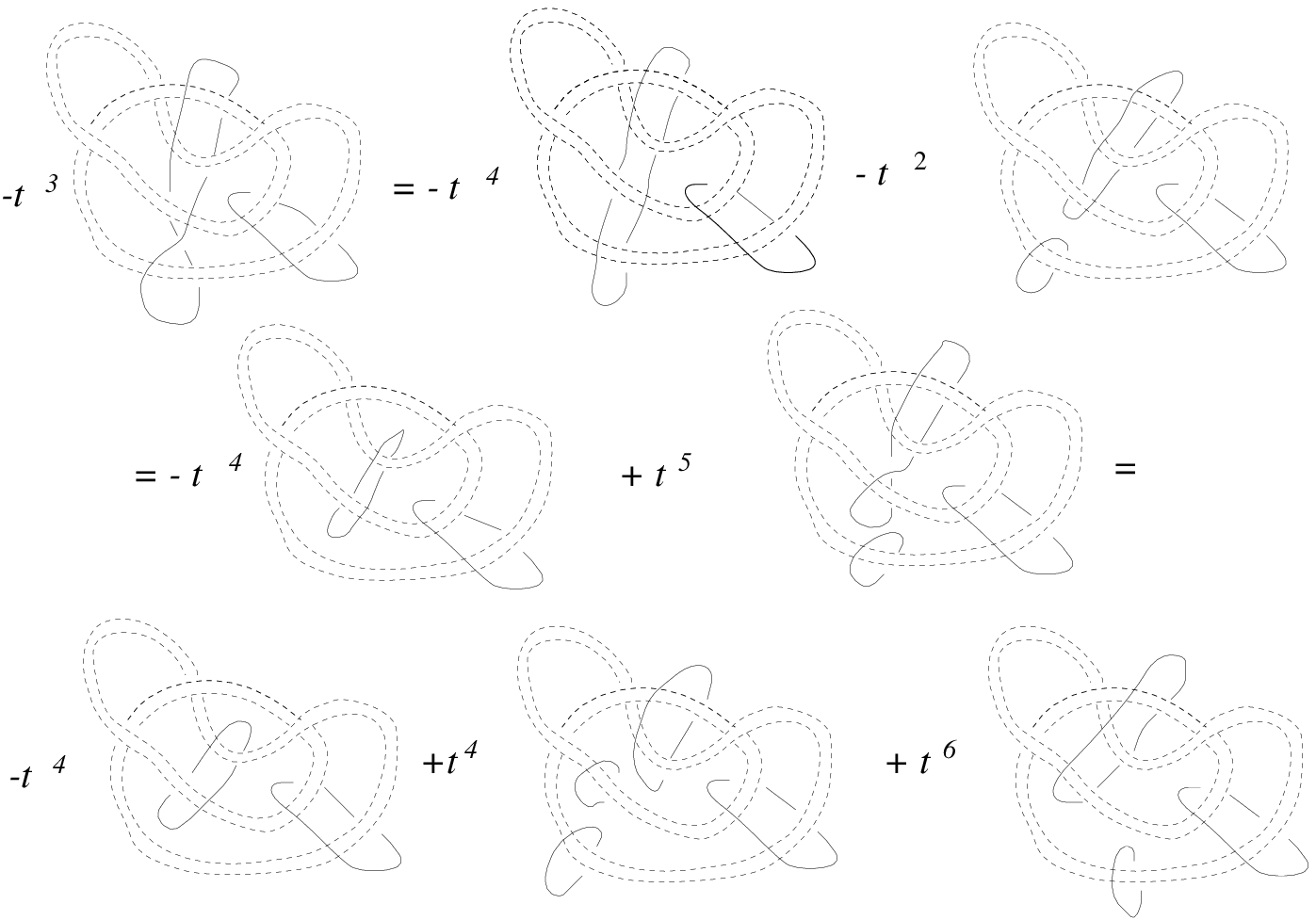}

Figure 5.  
\end{figure}

In the last line of Fig. 5, the second skein is equal to $t^4x^2y^2$ and
the third to $t^6x^2y$. Since we already computed $y^2$, it only
remains to compute the remaining one skein. 
To this end we write a Kauffman bracket skein relation as in
Fig. 6. The first diagram from the right side can be easily 
transformed like in  Fig. 7, and applying the skein relation we
get that this is equal to $-t^{-2}x^2-t^{-4}y$. 
For the diagram  on the left we use the computation from  Fig. 8 to conclude
 that is equal to $t^7+t^3-2t^3x^2-tx^2y$. Hence 
the missing skein is equal to $t^8+t^4-2t^4x^2+x^2+t^{-2}y-t^4x^2y$. 
In the end we deduce that $y^3=t^{10}+t^6-3t^6x^2+t^6x^4+2y+t^4y-
2t^4x^2y-t^4x^2y+t^4x^4y$, and by grouping terms we get the formula
from the statement.

\begin{figure}[htbp]
\centering
\leavevmode
\epsfxsize=3.5in
\epsfysize=1.2in
\epsfbox{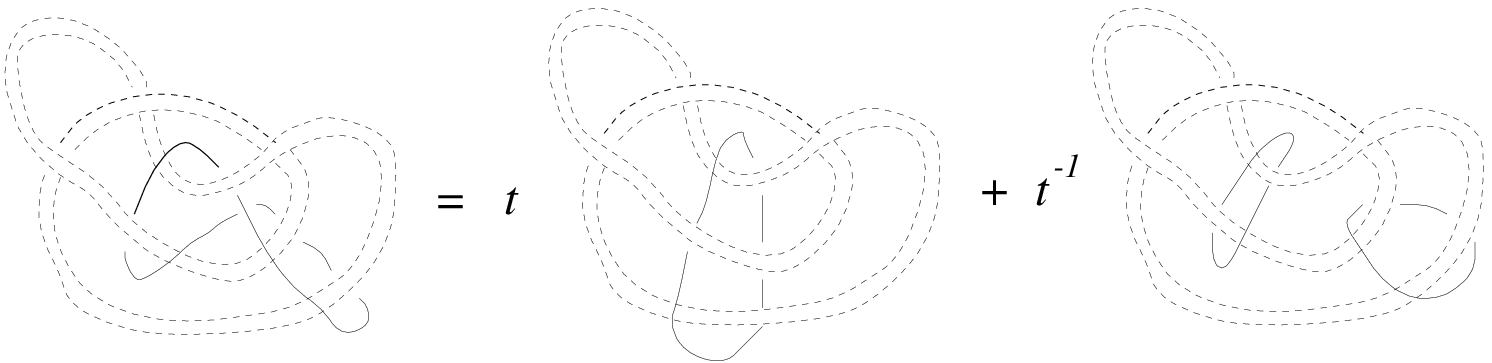}

Figure 6.  

\centering
\leavevmode
\epsfxsize=5in
\epsfysize=1.2in
\epsfbox{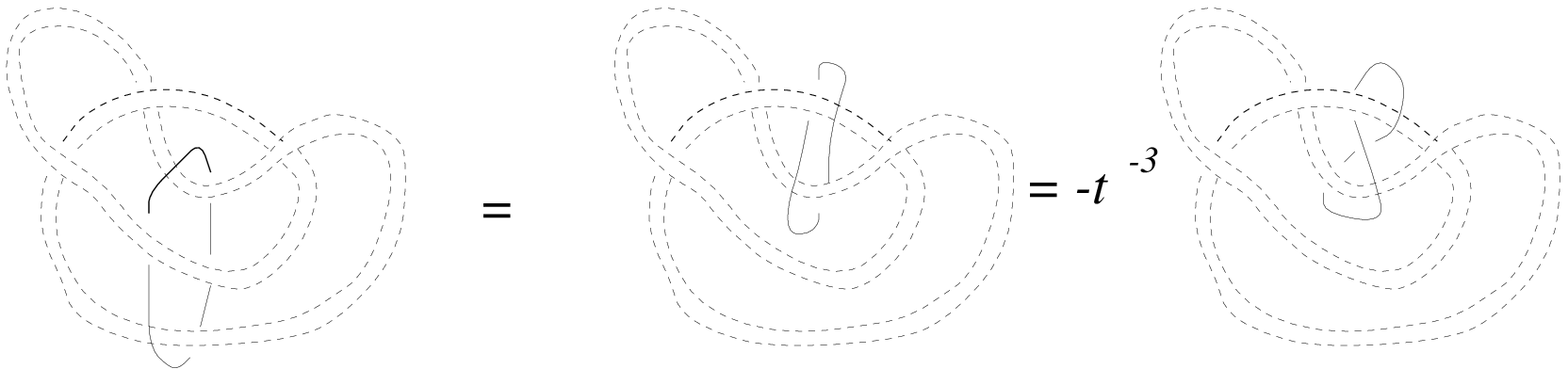}

Figure 7.  

\centering
\leavevmode
\epsfxsize=5in
\epsfysize=2.5in
\epsfbox{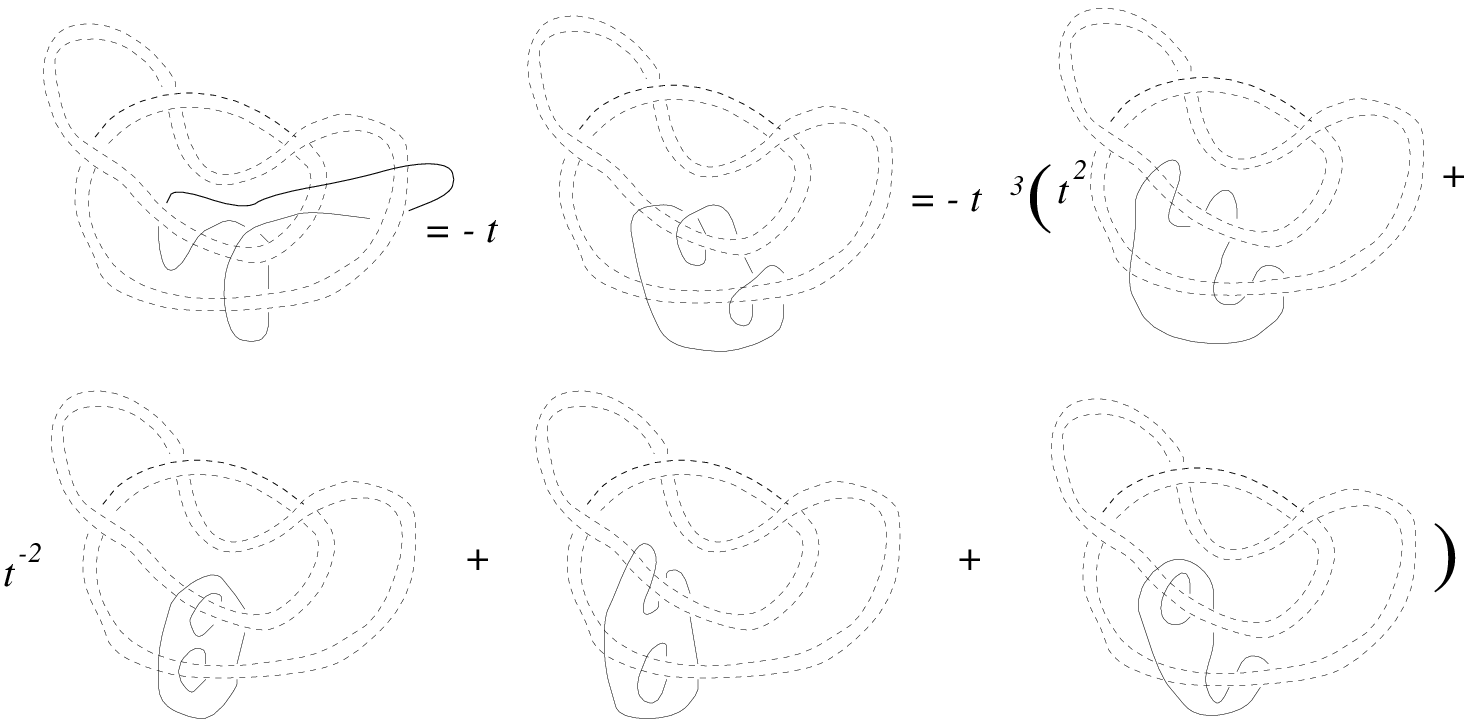}

Figure 8.  
\end{figure}

\qed

\begin{lemma}
The following equalities hold 
\begin{eqnarray*}
%& & a).\: \pi((1,-3)_T)=t^3S_3(x)+t^{-1}S_1(x)+(t+t^{-3})S_1(x)y.\\
%& & b). \: \pi((1,-2)_T)=t^4S_4(x)+S_0(x)+t^2S_2(x)y+t^{-2}S_0(x)y.
& & a).\: \pi((1,0)_T)=t^6S_6(x)-t^{2}S_0(x)+t^4S_4(x)y-S_0(x)y.\\
& & b). \: \pi((1,-1)_T)=t^5S_5(x)+t^3S_3(x)y.
\end{eqnarray*}
\end{lemma}

\proof

Using the Kauffman bracket skein relation we get the expansion
from Fig. 9. Using the computation from Fig. 7 we see that the first 
diagram is equal to $(-t^{-2}m^2-t^{-4}y)^3$. The second and the fourth
diagram are equal to $(-t^2-t^{-2})^2$ and $(-t^2-t^{-2})$, respectively.
It remains to find the value of the third diagram. The computations
from Fig. 10 and Fig. 11, combined with that of Fig. 7, show
 that the desired value of diagram is   
$1+t^{-4}+t^{-6}x^2y+t^{-8}x^4+t^{-10}y$.
On the other hand, the first skein from Fig. 7 is almost the 
image in the skein module of the knot complement of the $(1,0)$ curve
on the torus. To get the image of the curve  multiply
this skein by the framing factor of $-t^9$, induced
by the twistings of the $(1,0)$ curve around the torus, and a) follows.
The proof of b). is  analogous, even simpler, and is left to
the reader.
 
\begin{figure}[htbp]
\centering
\leavevmode
\epsfxsize=3.5in
\epsfysize=1.9in
\epsfbox{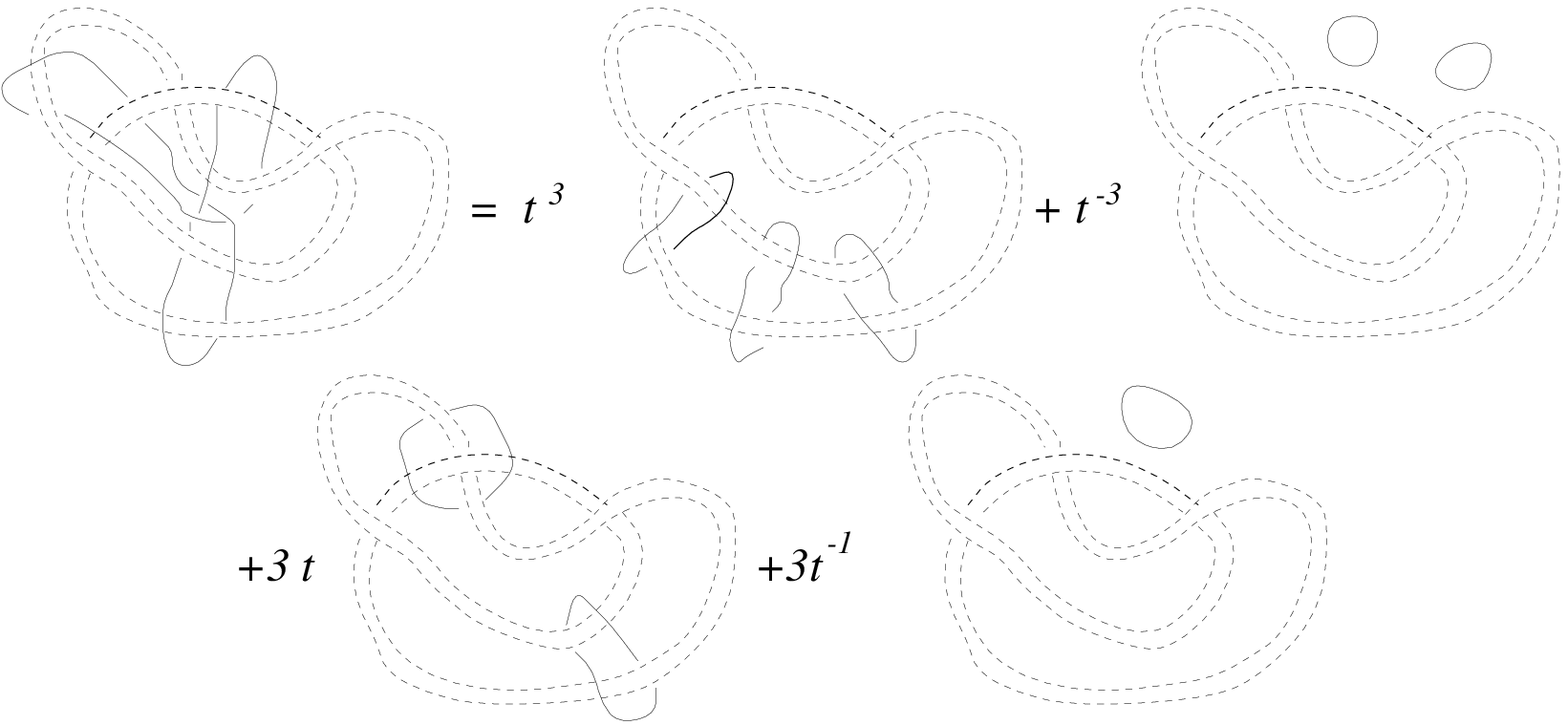}

Figure 9.  

\centering
\leavevmode
\epsfxsize=3.5in
\epsfysize=1.2in
\epsfbox{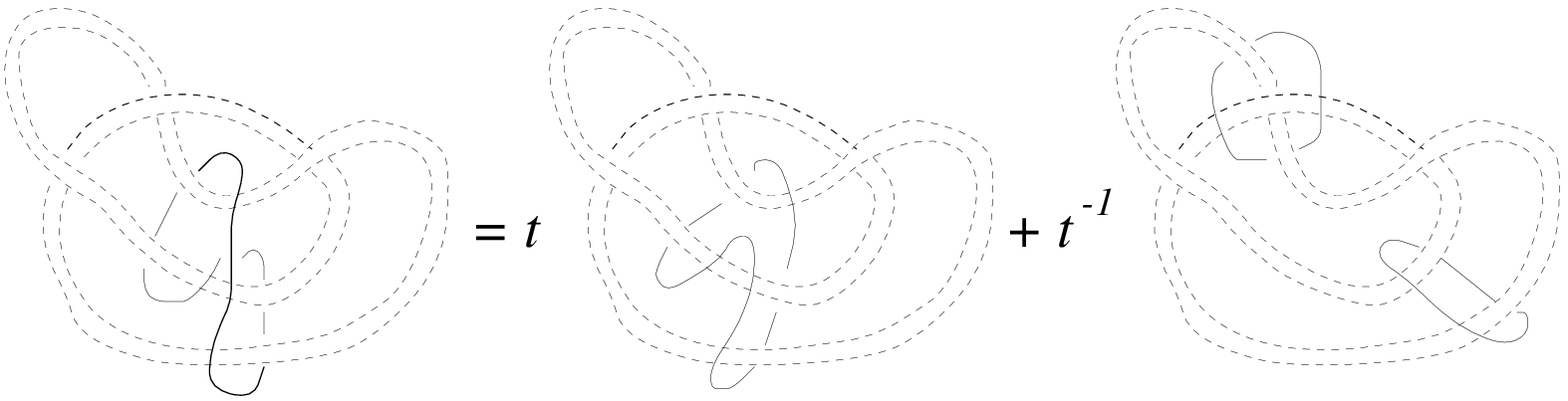}

Figure 10.  

\centering
\leavevmode
\epsfxsize=3.5in
\epsfysize=3.5in
\epsfbox{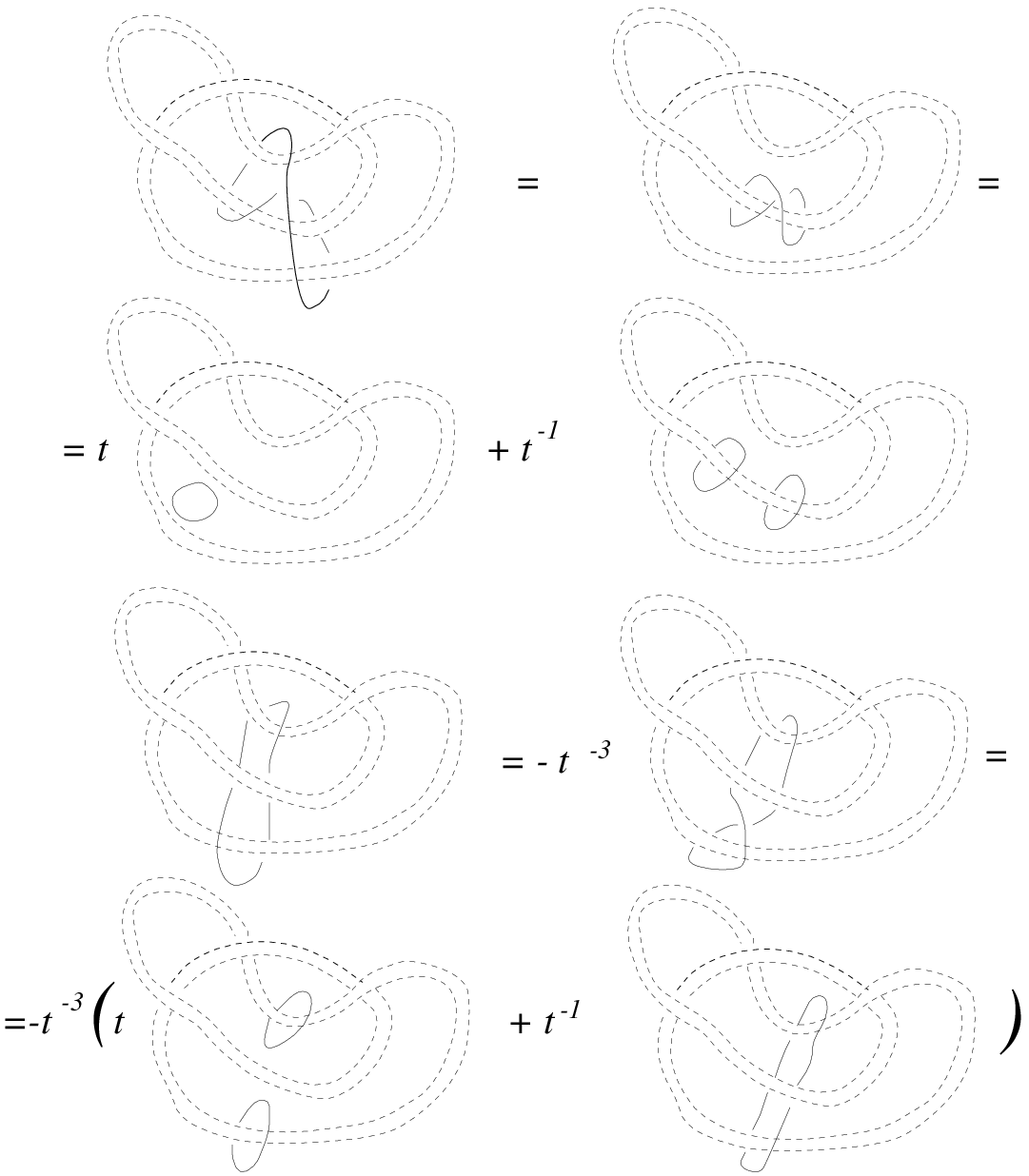}

Figure 11.  
\end{figure}

\qed

\begin{lemma}
For any $q\in {\mathbb Z}$ one has
\begin{eqnarray*}
\pi((1,q)_T)=t^{q+6}S_{q+6}(x)-t^{q+2}S_q(x)+
t^{q+4}S_{q+4}(x)y-t^qS_q(x)y.
\end{eqnarray*}
\end{lemma}

\proof
The proof is by induction on $q$, based on Lemma 3 and
the {\em product-to-sum} formula:
\begin{eqnarray*}
(1,q)_T*(0,1)_T=t(1,q+1)_T+t^{-1}(1,q-1)_T.
\end{eqnarray*}
Here we use the fact that $\pi((0,1)_T)=x$.  
\qed

An element of the skein module of a manifold is called 
peripheral if it is the image through the inclusion
map of an element of the skien algebra of the cylinder over the
boundary. A skein module is peripheral if all of its elements
are peripheral.
As a consequence of
Lemma 4 we get:

\begin{lemma}
In $K_t(M)$, one has  
\begin{eqnarray*}
(t^4-t^{-4})y=\pi \left(t^4(1,-4)_T-t^{-2}(1,-2)_T+
t^2(0,4)_T-t^6(0,2)_T-t^6+t^{-2}\right).
\end{eqnarray*}
So  if $t$ is not an eighth root of unity, then $y$ is peripheral,
hence the skein module $K_t(M)$ is peripheral.
\end{lemma}

Let $m,n$ be two integers. Denote
 by $x_{n,m}$ the skein in the mapping cylinder of a pair
of pants given in Fig. 12. Here if $m$ or $n$ is negative, the
kinks are in the other direction.
\begin{figure}[htbp]
\centering
\leavevmode
\epsfxsize=1.6in
\epsfysize=1.4in
\epsfbox{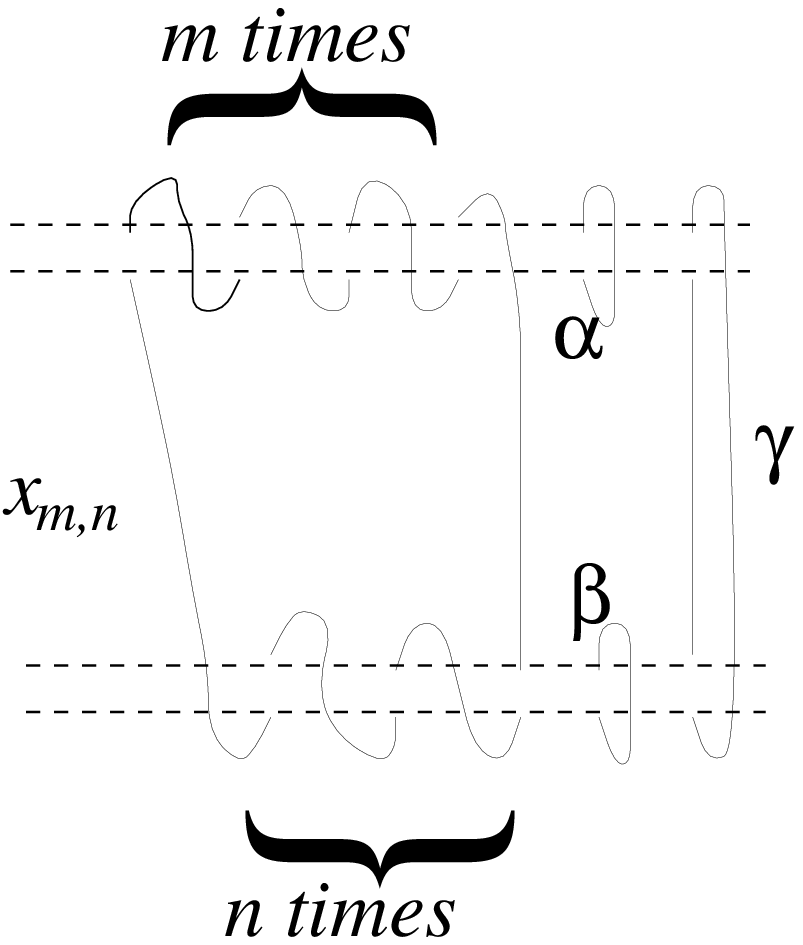}

Figure 12.  
\end{figure}
Recall that the mapping cylinder of a pair of pants is the free module
with basis $\alpha ^m\beta ^n\gamma ^p$, $m,n,p\geq 0$ where 
the curve $\alpha $, $\beta $ and $\gamma $ are shown in Fig. 12.

\begin{lemma}
\begin{eqnarray*}
x_{m,n} & = & (-t^{-2})^{m+n-1}S_m(\alpha)S_n(\beta)-
(t^{-2})^{m+n}S_{m-1}(\alpha)
S_{n-1}(\beta )\gamma +\\
& & (-t^{-2})^{n+m+1}S_{m-2}(\alpha )S_{n-2}(\beta).
\end{eqnarray*}
\end{lemma}

\proof
Sliding the strand one might produce a kink
which, when resolved via the Kauffman bracket skein relation,
produces two  configurations with less complexity.
Hence an induction on both $m$ and $n$ works to prove the formula.
\qed

\begin{lemma}
For $q\in {\mathbb Z}$ one has
\begin{eqnarray*}
(1,q)_T\cdot y=t^qS_{q-2}(x)-t^{q+8}S_{q+6}(x)+t^{q-2}S_{q-2}(x)y-
t^{q+6}S_{q+4}(x)y.
\end{eqnarray*}
\end{lemma}

\proof
 Note that we have the equality from Fig. 13. 
If we denote by $u_{n-3}$ the skein from Fig. 14 multiplied by 
$(-t^3)^n\cdot (-t^9)$, then the image of $(1,p)_T\cdot y$ in $K_t(M)$ is
obtained from the sum in Fig. 13 by capping each term of the sum
with $u_{p-3}$. Let us denote by $a_1,a_2,\cdots ,a_{12}$ be the 
twelve skeins of the sum respectively, and let us denote the
capping operation by $\circ $. Then one can easily
see that 
\begin{eqnarray*}
& & x_1\circ u_{n}=\pi ((1,n)_T)\\
& & x_2\circ u_n=(-t^{-3})x\pi ((1,n+1)_T)\\
& & x_3\circ u_n=(-t^3)x\pi ((1,n-1)_T).
\end{eqnarray*}
\begin{figure}[htbp]
\centering
\leavevmode
\epsfxsize=4.5in
\epsfysize=3.2in
\epsfbox{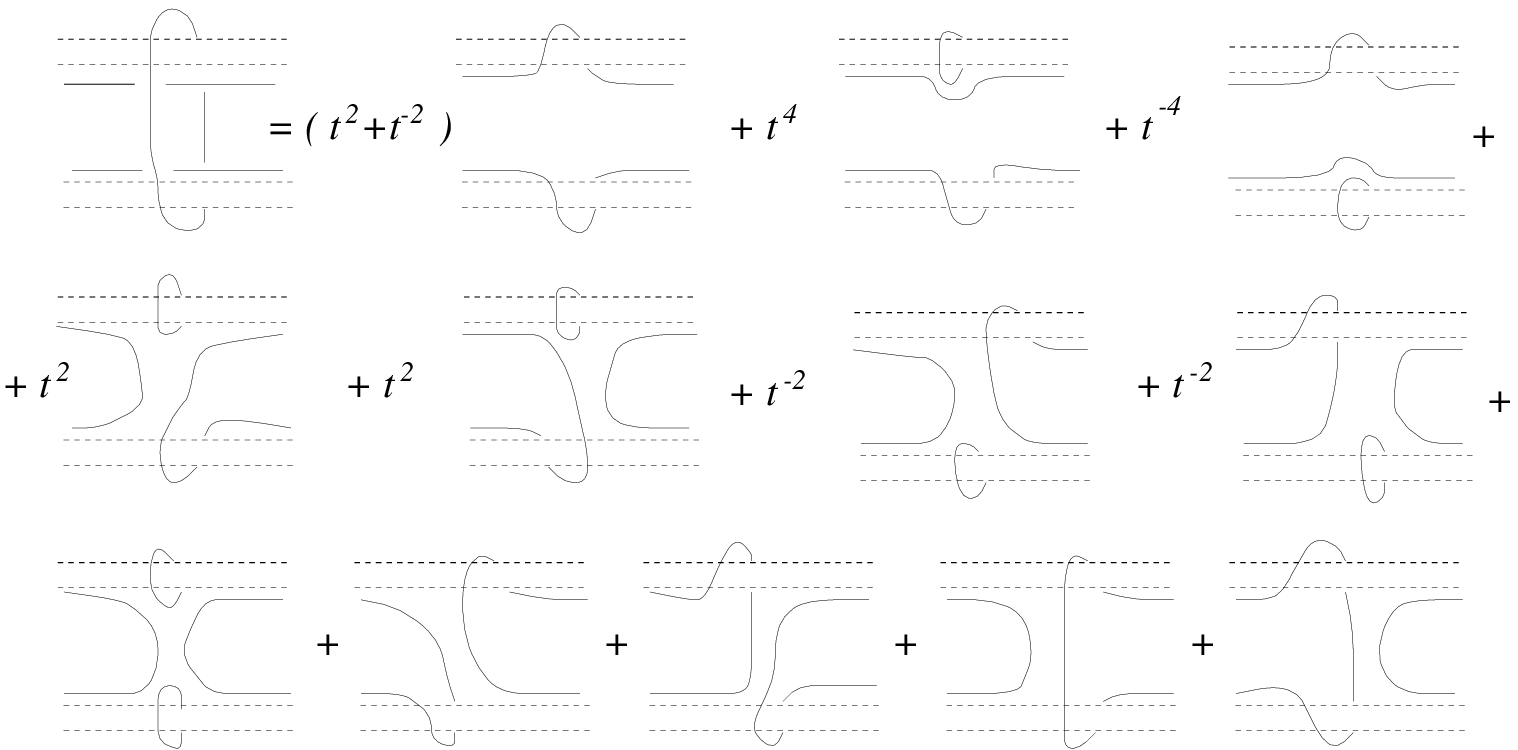}

Figure 13.  

\centering
\leavevmode
\epsfxsize=2.5in
\epsfysize=.9in
\epsfbox{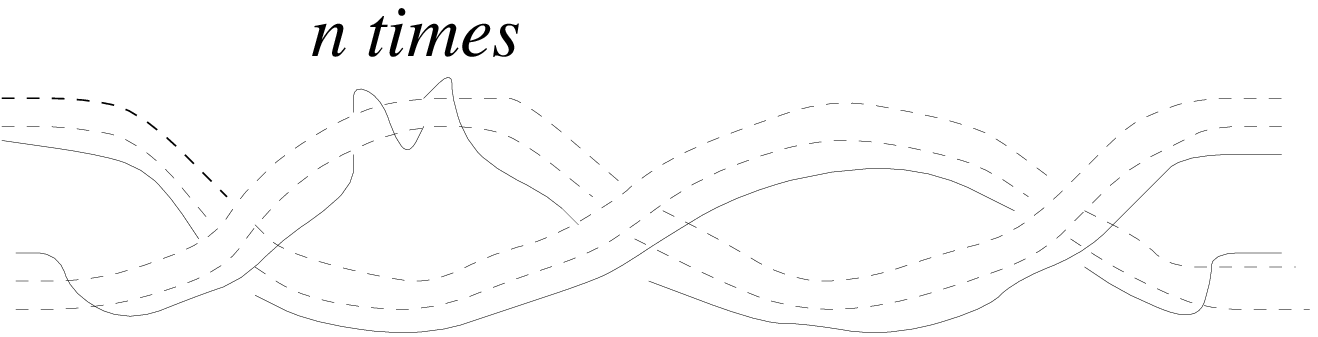}

Figure 14.  
\end{figure}
Also using Lemma  6 we get
\begin{eqnarray*}
& & x_4\circ u_n=t^{n+6}x^2S_{n+6}(x)+t^{n+4}xS_{n+5}(x)y\\
& & x_5\circ u_n=t^{n+6}xS_{2}(x)S_{n+5}(x)+t^{n+4}x^2S_{n+4}(x)+
t^{n+2}xS_{n+3}(x)\\
& & x_6\circ u_n= t^{n+10}xS_{n+5}(x)-t^{n+6}xS_{n+3}(x)\\
& & x_7\circ u_n= t^{n+10}x^2S_{n+4}(x)+t^{n+8}xS_{n+3}(x)y\\
& & x_8\circ u_n=-t^{n+8}x^3S_{n+5}(x)-t^{n+6}x^2S_{n+4}(x)y\\
& & x_9\circ u_n=-t^{n+8}xS_{n+5}(x)-t^{n+6}S_{n+4}(x)y\\
& & x_{10}\circ u_n=-t^{n+8}xS_{n+5}(x)-t^{n+6}S_{n+4}(x)y\\
& & x_{11}\circ u_n=-t^{n+8}S_{n+6}(x)+t^{n+4}S_{n+4}(x)\\
& & x_{12}\circ u_n=-t^{n+8}S_2(x)S_{n+4}(x)-t^{n+6}xS_{n+3}(x)y-
t^{n+4}S_{n+2}(x).
\end{eqnarray*}
On the other hand, $xS_n(x)=S_{n+1}(x)+S_{n-1}(x)$ and 
\begin{eqnarray*}
x\pi((p,q)_T)=\pi ((p,q)_T*(0,1)_T)=t^{p}\pi ((p, q+1)_T)+
t^{-p}\pi ((p,q-1)_T).
\end{eqnarray*}
Using these observations and Lemma 4, we see how to write all
right-hand sides in the above relations in the basis
of $K_t(M)$. A lengthy computation with many cancellations
produces the formula from the statement.
\qed

Define $\epsilon _k$ to be equal to $1$ if $k$ is odd, and $0$ if 
$k$ is even.

\begin{theorem}
If $p>0$ then
\begin{eqnarray*}
 & & \pi ((p,q)_T)=t^{6p^2+pq}S_{q+6p}(x)+t^{6p^2+pq-2}S_{q+6p-2}(x)y+\\
& & \sum_{k=1}^{2p-1}(-1)^{\left[\frac{k+1}{2}\right]}t^{6p^2-\left(\left[
\frac{k}{2}\right]+1\right)\left(6\left[\frac{k}{2}\right]+4\epsilon _k
\right)+pq}S_{q+6p-3k+\epsilon _k-2}(x)(1+t^{-2}y).
\end{eqnarray*}
The action of $K_t({\mathbb T}^2\times I)$ on $K_t(M)$ is 
determined by 
\begin{eqnarray*}
& & (p,q)_T\cdot y=-t^{6p^2+pq+2}S_{q+6p}(x)-t^{6p^2+pq}S_{q+6p-2}(x)y+\\
& & \sum_{k=1}^{2p-2}(-1)^{\left[\frac{k-1}{2}\right]}t^{6p^2-\left(\left[
\frac{k}{2}\right]+1\right)\left(6\left[\frac{k}{2}\right]+4\epsilon _k
\right)+pq}S_{q+6p-3k+\epsilon _k-2}(x)(t^2+y)+\\
& & (-1)^{p-1}S_{q}(x)t^{-2q+pq}(
t^2+y).
\end{eqnarray*}
\end{theorem}

\proof
The two equations are proved  simultaneously by induction. 
By Lemma 4 and Lemma 7, they are true for $p=1$, and the inductive
step is based on the {\em product-to-sum} formula
\begin{eqnarray*}
(1,0)_T*(p,q)_T=t^{p}(p+1,q)_T+t^{-p}(p-1,q)_T.
\end{eqnarray*}
and the fact that, by Lemma 7,  
\begin{eqnarray*}
(1,0)_T\cdot y=-1-t^8S_6(x)-t^{-2}y-t^6S_4(x)y.
\end{eqnarray*} 
\qed

\section{The noncommutative A-ideal of the left-handed trefoil}

Fix  $t$ is a complex number, which is
not an eighth root of unity. The considerations below hold also in
any coefficient ring in which $t^8-1$ is invertible.
For $m\geq 1$, denote by $I_t^m(K)$ the intersection of 
$I_t(K)$ with the linear span of $(p,q)_T$, $0\leq p\leq m$, 
$q\in {\mathbb Z}$. 

\begin{lemma}
Every element in $I_t^{1}(K)$ is of the form $p((0,1))*\tau $, where
\begin{eqnarray*}
\tau = (1,-5)_T-t^{-8}(1,-1)_T+t^{-3}(0,5)_T-t(0,1)_T.
\end{eqnarray*}
and $p$ is a polynomial with complex coefficients . 
\end{lemma}

\proof 
There is a part of $I_t^1(K)$ that arises via Lemmas $4$ and $5$, namely
that spanned by the elements
\begin{eqnarray*}
\phi_q & = & (t^4-t^{-4})((1,q)_T-t^{q+6}(0,q+6)_{JW}+t^{q+2}(0,q)_{JW})-\\
& & (t^{q+4}(0,q+4)_{JW}-t^q(0,q)_{JW})*(t^4(1,-4)_T-t^{-2}(1,-2)_T+\\
& & +t^2(0,4)_T
-t^6(0,2)_T-t^6+t^{-2}).
\end{eqnarray*}
For example, for $q=-5$ we get $\tau$, and for $q=-6$, we get the element in 
the kernel
\begin{eqnarray*}
& & (1,-6)_T+t^{-2}(1,-4)_T-t^{-8}(1,-2)_T-t^{-10}(1,0)_T+t^{-4}(0,6)_T+\\
& & t^{-4}(0,4)_T-(0,2)_T-2,
\end{eqnarray*}
and it is not hard to see that it is equal to $(0,1)_T*\tau$. 
On the other hand, using the {\em product-to-sum} formula, we can write
$(0,1)_T*\phi_q=t^{-1}\phi_{q+1}+t\phi_{q-1}$ so an induction
in both direction shows that all elements $\phi _q$ are of the form
$p((0,1)_T)*\tau $.
Note also that  if $q\geq 0$, then
the terms with extreme second coordinate
 that appear in the formula of $\phi _q$ are 
$t^{-2}(1,q+2)_T$ and $t^{2q+2}(1,-q-8)_T$. 

Among all elements in the part of $I_t^1(K)$ not spanned by $\phi _q$'s, 
choose one, $\psi$,  
such that the maximum of $|q+1|$ with  $(1,q)_T$ appearing
in the writing of $\psi$, is minimal.  
Let $a(1,q+2)_T +(b(1,-q-4)_T$ be the part where the maximum
is attained (here $q\geq 0$). Note that in $K_t(M)$,
 $\pi(\psi)$ had the coordinate
of $S_{q+6}(x)y$ equal to $at^{q+6}+bt^{q-8}$, hence $a$ and $b$ are
in the same proportion as those of $\phi _q$. But then by subtracting
a multiple of $q$, we can eliminate them, contradicting the minimality
of $\psi$. This shows that there are no other elements in
$I_t^1(K)$, and the lemma is proved.
\qed

\begin{theorem}
The ideal $I_t(K)$ is generated by 
\begin{eqnarray*}
& & (1,-5)_T-t^{-8}(1,-1)_T+t^{-3}(0,5)_T-t(0,1)_T,\\
& & (2,-6)_T-(t^6+t^{-6})(1,0)_T+(t^4+t^{-4})(1,-6)_T+(0,6)_T-2(t^4+t^{-4}),\\
& & (2,-7)_T-t^{-5}(1,-7)_T+(t^{-5}-t^{-1})(1,-3)_T-t^5(1,-1)_T+\\
& & + (t^2-t^{-2})(0,3)_T-t^{-6}(0,1)_T.
\end{eqnarray*}
\end{theorem}

\proof 
Let us denote by ${\mathcal J}$ the ideal generated by the three skeins
from the statement. The fact that ${\mathcal J}\subset I_t(K)$ follows from 
 from Lemmas 4, 5 and 8. Let us prove the equality. First, we check that
there are no other generators in $I_t^2(K)$. By Lemma 8, we only need
to check the part of $I_t^2(K)$ of elements that actually
contain $(2,q)_T$'s. But using the {\em product-to-sum} formula
\begin{eqnarray*}
(2,q)_T*(0,1)_T=t^{2}(2,q+1)_T+t^{-2}(2,q-1)_T
\end{eqnarray*} 
we can prove inductively that each such element in $I_{t}^2(K)$ can be 
reduced modulo ${\mathcal J}$ to one in $I_t^1(K)$, and the latter is
in ${\mathcal J}$ by Lemma 8. 

We will prove that $I_t^m(K)$ is contained in ${\mathcal J}$
by induction on $m$. 
First note that $I_t(K)$ is spanned by the elements that arise by
using the first formula from Theorem 1, in which we replace $x$ by
$(0,1)_T$ and $y$ by
its value given in Lemma 5.  Indeed, the right-hand side of the formulas
produces only elements in $(0,n)_T$, $(1,n)_T$, we can eliminate all
$(p,q)_T$'s that appear in the writing of an elements in the kernel
using them and reduce such an element modulo ${\mathcal J}$ to
an element that only contains $(p,q)_T$ with $p\leq 2$, and then use
the first part of the proof.

On the other hand, using the {\em product-to-sum} formula
\begin{eqnarray*}
(1,0)_t*(p,q)_T=t^{q}(p+1,q)_T+t^{-q}(p-1,q)_T
\end{eqnarray*}
we see that inductively we can reduce any element in $I_t^m(K)$ to
one in $I^2_t(K)$ modulo ${\mathcal J}$, and the conclusion follows.  
\qed

\begin{theorem}
The ideal ${\mathcal A}_t(K)$ is generated by 
\begin{eqnarray*}
& & [m^4(l+t^{10})-t^{-4}(l+t^{2})](l-t^6m^6),\\
& & (l+t^{24})(l+t^{10})(l+t^2)(l-t^6m^6),\\
& & (m^2-t^{-22})(l+t^{10})(l+t^2)(l-t^6m^6).
\end{eqnarray*}
\end{theorem}

\proof
The first generator of $I_t(K)$ given in Theorem 2 gives rise 
to the element 
\begin{eqnarray*}
 t^5(lm^{-5}+l^{-1}m^5)-t^{-7}(lm^{-1}+l^{-1}m)+t^{-3}(m^5+m^{-5})-t(m+m^{-1})
\end{eqnarray*}
in the subring of the noncommutative torus consisting 
of trigonometric polynomials. After multiplying by $t^{11}$ this 
element contracts to
\begin{eqnarray*}
t^{-4}l^2+t^{16}m^{10}-t^{-16}l^2m^4-t^{4}m^6+t^{-2}lm^{10}+t^{-2}l-t^{2}lm^6
-t^{2}lm^4
\end{eqnarray*}
 in the quantum plane. It is not hard to see that it
factors as 
$[m^4(l+t^{10})-t^{-4}(l+t^{2})](l-t^6m^6)$.
The other two generators of $I_t(K)$ give rise to two
more generators of ${\mathcal A}_t(K)$. The ones from the statement
are obtained from these two after some algebraic manipulations
which involve also the first generator.  
\qed

Note the presence of the factor $(l-t^6m^6)$, which is the noncommutative
analogue of the factor 
of the classical A-polynomial
that stands for the irreducible $SL(2,C)$-representations
of the fundamental group.

Now let  $t=-1$. Note  the presence of a discontinuity, 
due to Lemma 5. In this case, the same arguments apply {\em mutatis
mutandis} to prove the following two results.

\begin{theorem}
The ideal $I_{-1}(K)$ is generated by $(1,-4)_T-(1,-2)_T+(0,4)_T-(0,2)_T-2$
and $(2,-6)_T-(0,6)_T$.
\end{theorem}

\begin{theorem}
The ideal ${\mathcal A}_{-1}(K)$ is generated by
$(l^2-1)(l+1)(l-m^6)$ and $(m^2-1)(l+1)(l-m^6)$.  
\end{theorem}

Note that the classical A-polynomial is obtained by replacing
$l$ by $-l$ and $m$ by $-m$, and then taking the common factors 
of the two generators (i.e. eliminating the embedded  primes).
The change of sign in the variables is due to the fact that
the relationship between skein modules and character varieties
is established by the negative of the trace.

\section{The case of the right-handed trefoil}

When taking the mirror image the $(p,q)$ curve in the boundary
of the left-handed trefoil becomes the $(p,-q)$ curve on the boundary
of the right-handed trefoil. Also, in the Kauffman 
bracket skein relation, taking the mirror image changes $t$ to $t^{-1}$. 
So for $K$ the right-handed trefoil knot we get the following results.

\begin{theorem}
If $p>0$ then
\begin{eqnarray*}
 & & \pi ((p,q)_T)=t^{-6p^2+pq}S_{-q+6p}(x)+t^{-6p^2+pq+2}S_{-q+6p-2}(x)y+\\
& & \sum_{k=1}^{2p-1}(-1)^{\left[\frac{k+1}{2}\right]}t^{-6p^2+\left(\left[
\frac{k}{2}\right]+1\right)\left(6\left[\frac{k}{2}\right]+4\epsilon _k
\right)+pq}S_{-q+6p-3k+\epsilon _k-2}(x)(1+t^{2}y).
\end{eqnarray*}
The action of $K_t({\mathbb T}^2\times I)$ on $K_t(M)$ is 
determined by 
\begin{eqnarray*}
& & (p,q)_T\cdot y=-t^{-6p^2+pq+2}S_{-q+6p}(x)-t^{-6p^2+pq}S_{-q+6p-2}(x)y+\\
& & \sum_{k=1}^{2p-2}(-1)^{\left[\frac{k-1}{2}\right]}t^{-6p^2+\left(\left[
\frac{k}{2}\right]+1\right)\left(6\left[\frac{k}{2}\right]+4\epsilon _k
\right)+pq}S_{-q+6p-3k+\epsilon _k-2}(x)(t^{-2}2+y)+\\
& & (-1)^{p-1}S_{-q}(x)t^{-2q+pq}(
t^{-2}+y).
\end{eqnarray*}
\end{theorem}

\begin{theorem}
If $t$ is not an eighth root of unity, the ideal $I_t(K)$ is generated by 
\begin{eqnarray*}
& & (1,5)_T-t^{8}(1,1)_T+t^{3}(0,5)_T-t^{-1}(0,1)_T,\\
& & (2,6)_T-(t^6+t^{-6})(1,0)_T+(t^4+t^{-4})(1,6)_T+(0,6)_T-2(t^4+t^{-4}),\\
& & (2,7)_T-t^{5}(1,7)_T+(t^{5}-t)(1,3)_T-t^{-5}(1,1)_T-\\
& & - (t^2-t^{-2})(0,3)_T-t^{6}(0,1)_T.
\end{eqnarray*}
\end{theorem}

\begin{theorem}
If $t$ is not an eighth root of unity, 
the ideal ${\mathcal A}_t(K)$ is generated by 
\begin{eqnarray*}
& & [m^4(l+t^{10})-t^{-4}(l+t^{2})](lm^6-t^6),\\
& & (l+t^{24})(l+t^{10})(l+t^2)(lm^6-t^6),\\
& & (m^2-t^{-22})(l+t^{10})(l+t^2)(lm^6-t^6).
\end{eqnarray*}
\end{theorem}

\begin{theorem}
The ideal $I_{-1}(K)$ is generated by $(1,4)_T-(1,2)_T+(0,4)_T-(0,2)_T-2$
and $(2,6)_T-(0,6)_T$,
and the ideal ${\mathcal A}_{-1}$ is generated by
$(l^2-1)(l+1)(lm^6-1)$ and $(m^2-1)(l+1)(lm^6-1)$  
\end{theorem}

\end{document}